\documentclass[12pt]{article}
\usepackage{amsmath,amsthm}
\usepackage{amssymb}
\usepackage{enumerate}
\usepackage{hyperref}
\newtheorem{theorem}{Theorem}[section]

\newtheorem{lemma}[theorem]{Lemma}
\newtheorem{remark}[theorem]{Remark}
\def\proofof{\goodbreak\noindent{\sc Proof of }\nobreak}
\def\endproof{\par\nobreak\hbox to \hsize{\hfil\vrule width 5pt height 5pt}\goodbreak\vskip 3pt}
\def\endproofof{\par\nobreak\hbox to \hsize{\hfil\vrule width 5pt height 5pt}\goodbreak\vskip 3pt}
\def\bbb{\mathbb}
\def\bR{{\bbb R}}
\def\cB{{\cal B}}
\def\cD{{\cal D}}
\def\cE{{\cal E}}
\newcommand{\sptext}[3]{\hspace{#1 em}\mbox{#2}\hspace{#3 em}}
\newcommand{\equa}{\begin{eqnarray*}}
\newcommand{\tion}{\end{eqnarray*}}
\newcommand{\vare}{\varepsilon}
\def\y{\lbrack\!\lbrack} 
\def\yy{\rbrack\!\rbrack}
\def\sm{\setminus}

\allowdisplaybreaks
\begin{document}

\title{Haar type and Carleson Constants}
\author{Stefan Geiss \and
        Paul F.X. M\"uller
\thanks{Research of both authors supported in part by FWF Pr. Nr. P150907-N01.}}

\maketitle
\begin{abstract} 
For a collection $\cE $ of dyadic intervals,
a Banach space $X$, and
$p\in (1,2]$ we assume the 
upper $\ell^p$ estimates
\[      \left\| \sum _{I \in \cE} x_I h_I/|I|^{1/p} \right\|^p_{L^p_X} 
    \le c^p \sum _{I \in \cE} \| x_I\|_X^p ,\]
where  $ x_I \in X $ and $h_I$ denotes the $L^\infty$ normalized Haar 
function supported on $I$. We determine the minimal requirement
on the size of $\cE $ so that these estimates imply that
$ X$ is of Haar type $p. $ The characterization is given in terms
of the Carleson constant of $\cE . $
\end{abstract}




\section{Introduction}

Let $X$ be a Banach space.
We fix a non-empty collection of dyadic intervals $\cE$
and assume the upper $\ell^p$ estimates 
\begin{equation}
\label{upper}  
        \left\| \sum _{I \in \cE} x_I \frac{h_I}{|I|^{1/p}} \right\|_{L^p_X} 
    \le c \left ( \sum _{I \in \cE} \| x_I\|_X^p \right )^\frac{1}{p}
\end{equation}
for finitely supported $(x_I)_{I\in\cE}\subset X$ and some $p\in (1,2]$, where
$h_I$ is the $L^\infty$ normalized Haar function supported on $I$.
The consequences for $X$, one may draw from \eqref{upper}, depend 
on the size and structure of the collection $\cE . $ 
For instance, if $\cE $ is a collection of pairwise disjoint dyadic 
intervals, then any Banach space  satisfies \eqref{upper}, hence it 
does not impose any restriction on $X.$
If, on the other hand,  we choose $\cE$ to be  the collection of {\it all}
dyadic intervals, then the upper $\ell ^p $ estimates \eqref{upper} 
simply state that $X$ is of Haar type $p;$ due to important  work of G. Pisier \cite{pis} 
the latter condition is equivalent to certain renorming properties of 
the Banach space $X$.
  
In this paper we ask how massive a collection 
$\cE$ has to be so that \eqref{upper}  implies that 
$X$ is of Haar type $p$. We give the  answer  
to  this question in terms of  the Carleson 
constant
defined by 
\begin{equation}\label{eqn:Carleson_constant}
\y \cE \yy := \sup_{I\in \cE} \frac{1}{|I|} \sum_
{ J \in \cE ,\, J \subseteq I}
|J|.
\end{equation}
The proof is based on the following well-known dichotomy:
either $\cE $ can be decomposed into finitely many collections 
consisting of ``almost disjoint'' dyadic intervals
or   $\cE $ contains large and densely packed blocks of dyadic intervals,
with arbitrary high degree of condensation. 
%

Initially we encountered the problem treated here in connection with
our efforts to obtain a vector valued version of E. M. Semenov's 
characterization of bounded operators rearranging the 
Haar system. See \cite{sem} and \cite{gm}.

\section{Preliminaries}

In the following we equip the unit interval $[0,1)$ with the Lebesgue measure
denoted by $|\cdot|$.
Let  $\cD $ denote the collection of dyadic intervals in $[0,1)$, i.e.
$I \in \cD$ provided that there exist $m\ge 0$ and $1\le k\le 2^m$ such that 
\[ I = [(k-1)/2^m , k/2^m), \]
and let 
\[ \cD _n := \{ I \in \cD :  |I| \ge 2^{-n} \}
   \sptext{1}{where}{1}
   n\ge 0. \]
The $L^\infty$ normalized Haar function supported on $I\in\cD$
is denoted by $h_I$, i.e. $h_I = -1$ in the right half of $I$ and 
$h_I = 1$  on the left half of $I .$
By $L_X^p$, $p\in [1,\infty)$, we denote the space of Radon random 
variables $f:[0,1)\to X$ such that
\[ \|f\|_{L_X^p } = \left(\int_0^1\|f(t)\|_X ^p dt \right)^{1/p} < \infty. \]

\paragraph{Haar type.}
Given $p\in (1,2]$, a Banach space $X$ is of Haar type $p$
provided that there exists a constant 
$c > 0 $ such that 
\[      \left\| \sum _{I \in \cD} x_I \frac{h_I}{|I|^{1/p}} \right\|_{L^p_X} 
    \le c \left ( \sum _{I \in \cD} \| x_I\|_X^p\right )^\frac{1}{p} \]
for all finitely supported families 
$(x_I)_{I \in \cD}\subset X$. We let $HT_p(X)$ be the infimum
of all possible $c>0$ as above.
The central result concerning Haar type is due to G. Pisier \cite{pis} and 
asserts that Haar type $p$ is equivalent to the fact that $X$ can be equivalently 
renormed such that the new norm has a modulus of smoothness of power type 
$p.$ For additional information see \cite{Deville_Godefroy_Zizler,piwe} and the 
references therein.

\paragraph{Carleson Constants.}
Let $\cE \subseteq \cD$ be a non-empty collection of dyadic 
intervals. The Carleson constant of $\cE$  is given by
equation (\ref{eqn:Carleson_constant}).
Next we define consecutive generations of $\cE $ and, using 
$ \y \cE \yy$, describe a dichotomy for $\cE $ known as 
the almost disjointification and condensation properties.
 
We define $G_0(\cE)$ to be the maximal dyadic intervals
of $\cE$ where {\em maximal} refers to inclusion.
Suppose, we have already defined 
$G_0(\cE )$,..., $G_p(\cE )$, then we form
\[ G_{p+1} (\cE) := G_0( \cE\sm  (G_0(\cE )\cup\dots \cup G_p(\cE ))  ). \]
Given $I\in \cD$, let $  I\cap \cE =\{ J \in \cE ,\, J \subseteq I\}$
and put
\[ G_k(I,\cE) := G_k(I\cap \cE)
   \sptext{1}{for}{1}
   k\ge 1. \]
Assume that $\y \cE \yy < \infty$ and that $M$ is the largest integer
smaller than $4\y \cE \yy+1$. Then
\begin{equation}\label{eqn:decomposition_E}
\cE_i := \bigcup_{k=0}^\infty G_{Mk + i} (\cE ) ,\quad 0 \le i \le M-1,
\end{equation}
consists of almost disjoint dyadic intervals, in the sense that 
for  $ I \in \cE_i ,$ 
\[ \sum_{J  \in G_1(I,\cE_i)}  |J| \le \frac {|I|}{2} 
\quad\text{and}\quad 
\sum _{ K \in I \cap \cE_i} |K| \le 2 |I| .
\]
Conversely, if $\y \cE \yy = \infty$, then for all $n \ge 1$ and $\vare\in(0,1)$
there exists a $K_0 \in \cE $ that is densely packed in the sense that  
\[ \sum_{J \in G_n ( K_0 , \cE )} |J| \ge (1- \vare) |K_0| .\]
Based on this we show in 
Lemma \ref{lemma:sufficient_condition_type}
that for any $n  $ the span $(h_I)_{I\in\cE}$
contains a system of functions,  with the same joint  distribution 
as the first $2^n$ elements of the Haar basis.

The proof of this basic dichotomy and some of its applications can be found 
in \cite[Chapter 3]{pfxm1}.


\section{Haar Type and Carleson Constants}
The main results of this note are Theorems \ref{main} and \ref{main_inverse}.
Combined they give an answer to the question as to which families 
of dyadic intervals $\cE $ will detect whether a Banach space $X$ has 
Haar type $p .$ The answer is a dichotomy: either $\y\cE \yy  <\infty ,$
then $ \cE $ does not detect any Haar type of any Banach space, 
or   $\y\cE \yy  =\infty ,$ then  $ \cE $ determines the Haar type 
$p$ constant exactly, for any Banach space $X$ and each $ 1 < p \le 2 . $
\smallskip

\begin{theorem}\label{main} 
Let $p\in (1,2]$ and $\cE\subseteq \cD$ be a non-empty collection of dyadic 
intervals. Then the following statements are equivalent:
\begin{enumerate}[$(1)$]
\item 
A Banach space $X$ is Haar type $p$ if there exists 
      $c>0$ such that for all finitely supported families 
      $(x_I)_{I\in\cD} \subset X$ one has
      \[     \left\| \sum _{I \in \cE} x_I \frac{h_I}{|I|^{1/p}} \right\|_{L^p_X} 
         \le c \left (\sum _{I \in \cE} \| x_I\|_X^p\right )^\frac{1}{p} ;\]
      the infimum over all such $c> 0 $ coincides with $HT_p(X).$
\item $ \y \cE \yy = \infty. $
\end{enumerate} 
\end{theorem}
\medskip

\begin{theorem}\label{main_inverse} 
Let $p\in (1,2]$, $\cE\subseteq \cD$ be a non-empty collection of dyadic 
intervals, and $X$ be a Banach space which is not of Haar type $p$. 
Then the following statements are equivalent:
\begin{enumerate}[$(1)$]
\item There exists a constant $c>0$ such that for all finitely supported 
      families $(x_I)_{I\in\cD} \subset X$ one has
      \[    \left\| \sum _{I \in \cE} x_I \frac{h_I}{|I|^{1/p}} \right\|_{L^p_X} 
         \le c \left (\sum _{I \in \cE} \| x_I\|_X^p\right )^\frac{1}{p}. \]
\item $ \y \cE \yy < \infty. $
\end{enumerate} 
\end{theorem}

Theorem~\ref{main} and Theorem~\ref{main_inverse} follow immediately 
from the following two lemmas (and the obvious fact that there are
Banach spaces without Haar type $p$ if $p\in (1,2]$).

\begin{lemma} 
Let $p\in (1,\infty)$, $\y \cE \yy < \infty$, and $X$ be a Banach space.
Then there is a constant $c_p>0$, depending at most on $p$, such that
one has
\[     \left\| \sum _{I \in \cE} x_I \frac{h_I}{|I|^{1/p}} \right\|_{L^p_X} 
   \le c_p \y \cE \yy^{1-\frac{1}{p}}
       \left (\sum _{I \in \cE} || x_I||_X^p\right )^{\frac{1}{p}} \]
for all  finitely supported $(x_I)_{I\in\cE} \subset X$.
\end{lemma}
\proof
Using (\ref{eqn:decomposition_E}), we write  
$ \cE = \cE_0 \cup \dots \cup \cE_{M-1} $ with  
$ M < 4 \y \cE \yy+1$ such that the collections $\{ A_I : I \in \cE_i \} $ 
of pairwise disjoint and measurable sets defined by 
\[ A_I := I \sm \bigcup_{ J \in G_1(I,  \cE_i) }  J , \quad I \in \cE_i, \]
satisfy 
\[ \frac{1}{2} |I| \le |A_I| \le |I|. \]
Since
\begin{eqnarray*}
      \left\| \sum_{  I \in \cE } x_I \frac{h_I}{|I|^{1/p}}\right\|_{L_X^p}  
&\le& \sum _{i = 0}^{M-1} \left \| \sum_{  I \in \cE_i } x_I \frac{h_I}{|I|^{1/p}}
      \right \|_{L_X^p} \\
&\le& M^{1-\frac{1}{p}} \left(  \sum _{i = 0}^{M-1} \left \|
      \sum_{  I \in \cE_i } x_I \frac{h_I}{|I|^{1/p}}
      \right \|_{L_X^p}^p \right)^{1/p}
\end{eqnarray*}
it is sufficient to prove 
\[     \left\| \sum_{  I \in \cE_i } x_I \frac{h_I}{|I|^{1/p}}\right\|_{L_X^p}
   \le c_p  \left ( \sum_{  I \in \cE_i } \|x_I\|_X^p \right )^\frac{1}{p} \quad \text{for} \quad  i\le M .\] 
But here we get that
\equa
&   & \left \| \sum_{I\in\cE_i} x_I \frac{h_I}{|I|^\frac{1}{p}} \right\|_{L_X^p}\\
& = & \left ( \sum_{K\in\cE_i} \int_{A_K}
      \left \| \sum_{I\in\cE_i} x_I \frac{h_I(t)}{|I|^\frac{1}{p}} 
      \right \|_X^p dt \right )^\frac{1}{p}  \\
& = & \left ( \sum_{K\in\cE_i} \frac{|A_K|}{|K|}\int_{A_K}
      \left \| \sum_{I\in\cE_i} x_I \left (\frac{|K|}{|I|}\right )^\frac{1}{p} h_I(t) 
      \right \|_X^p \frac{dt}{|A_K|} \right )^\frac{1}{p}  \\
&\le& \left ( \sum_{K\in\cE_i} \int_{A_K}
      \left \| \sum_{I\in\cE_i} x_I \left (\frac{|K|}{|I|}\right )^\frac{1}{p} h_I(t) 
      \right \|_X^p \frac{dt}{|A_K|} \right )^\frac{1}{p} \\
& = & \left ( \sum_{K\in\cE_i} \int_{A_K}
      \left \| \sum_{K\subseteq I\in\cE_i} x_I \left (\frac{|K|}{|I|}\right )^\frac{1}{p} h_I(t) 
      \right \|_X^p \frac{dt}{|A_K|} \right )^\frac{1}{p} \\
& = & \left ( \sum_{K\in\cE_i} \int_{A_K}
      \left \| \sum_{l=0}^{n(K)} x_{G_{-l}(K,\cE_i)}
               \left (\frac{|K|}{|G_{-l}(K,\cE_i)|}\right )^\frac{1}{p} h_{G_{-l}(K,\cE_i)}(t) 
      \right \|_X^p \frac{dt}{|A_K|} \right )^\frac{1}{p}\\
& = & \left ( \sum_{K\in\cE_i} \int_{A_K}
      \left \| \sum_{l=0}^\infty x_{G_{-l}(K,\cE_i)} \chi_{\{l\le n(K)\}}
               \left (\frac{|K|}{|G_{-l}(K,\cE_i)|}\right )^\frac{1}{p} h_{G_{-l}(K,\cE_i)}(t) 
      \right \|_X^p 
       \right . \\
&    & \left . \hspace*{26em}
       \frac{dt}{|A_K|} \right )^\frac{1}{p}
\tion
where $G_{-l}(K,\cE_i)$ form the maximal sequence of dyadic
intervals from $\cE_i$ such that
\[       K=G_0   (K,\cE_i) 
   \subset G_{-1}(K,\cE_i)
   \cdots
   \subset G_{-n(K)}(K,\cE_i) \]
and $G_{-n(K)}(K,\cE_i)$ is the unique maximal interval in $\cE_i$ containing $K$. Next we obtain an  upper estimate for the  last expression as follows:
\equa
&   & \hspace*{-2.5em}
      \sum_{l=0}^\infty \left ( \sum_{K\in\cE_i} \int_{A_K}
      \left \| x_{G_{-l}(K,\cE_i)} \chi_{\{l\le n(K)\}}
               \left (\frac{|K|}{|G_{-l}(K,\cE_i)|}\right )^\frac{1}{p} h_{G_{-l}(K,\cE_i)}(t) 
      \right \|_X^p \frac{dt}{|A_K|} \right )^\frac{1}{p} \\
& = & \sum_{l=0}^\infty \left ( \sum_{K\in\cE_i}
      \left \| x_{G_{-l}(K,\cE_i)} \chi_{\{l\le n(K)\}} \right \|_X^p 
      \frac{|K|}{|G_{-l}(K,\cE_i)|}
      \right )^\frac{1}{p}\\
& = & \sum_{l=0}^\infty \left ( \sum_{I,K\in\cE_i \atop G_{-l}(K,\cE_i)=I}
      \left \| x_I \right \|_X^p 
      \frac{|K|}{|I|}
      \right )^\frac{1}{p}\\
& = & \sum_{l=0}^\infty \left ( \sum_{I\in\cE_i}
      \left \| x_I \right \|_X^p 
      \frac{\sum_{K\in\cE_i \atop G_{-l}(K,\cE_i)=I} |K|}{|I|}
      \right )^\frac{1}{p}\\
&\le& \sum_{l=0}^\infty \left ( \sum_{I\in\cE_i}
      \left \| x_I \right \|_X^p 2^{-l} 
      \right )^\frac{1}{p}\\
& = & \left ( \sum_{l=0}^\infty 2^{-\frac{l}{p}} \right ) 
      \left ( \sum_{I\in\cE_i} \left \| x_I \right \|_X^p 
      \right )^\frac{1}{p}.
\tion
\endproof

Next we turn to the case $ \y \cE \yy = \infty$ for
which it is known that the Gamlen-Gaudet construction yields an approximation 
of the Haar system by appropriate 'blocks' of $(h_I)_{I\in \cE}$. 
The next lemma demonstrates that this construction fits perfectly  with our
Haar type inequalities. We carefully  avoid using  
 the unconditionality of the Haar system (and therefore the UMD property 
of Banach spaces) and  exhibit a system of functions 
with exactly the same joint distribution as  the ususal Haar basis, 
rather than to 
allow that the measures of the support of a Haar function and its 
approximation are related by uniformly bounded multiplicative constants.

\begin{lemma}\label{lemma:sufficient_condition_type}
Let $X$ be a Banach space, $p\in (1,2]$,
and $\cE $ be a collection of dyadic intervals such that 
\[ \y \cE \yy = \infty. \]
If there is a constant $c>0$ such that
\begin{equation}\label{eqn:E_type}
       \left\| \sum _{I \in \cE} x_I \frac{h_I}{|I|^{1/p}} \right\|_{L^p_X}
   \le c \left ( \sum _{I \in \cE} || x_I||_X^p \right )^\frac{1}{p}
\end{equation}
for all finitely supported families $(x_I)_{I\in\cE} \subset X$,
then $X$ is of Haar type $p$ with $HT_p(X)\le c$.
\end{lemma}

\begin{remark}
In Lemma \ref{lemma:sufficient_condition_type} the range $p\in (2,\infty)$
does not make sense since already $X=\bR$ does not have 
Rademacher type $p\in (2,\infty)$ and henceforth Haar type $p\in (2,\infty)$.
In other words, for $\y \cE \yy = \infty$ and $p\in (2,\infty)$
the inequality {\rm (\ref{eqn:E_type})} fails to be true.
\end{remark}

\proofof{Lemma  \ref{lemma:sufficient_condition_type}.}
Let $n \ge 1,$  $\delta\in (0,1)$, and
$\vare= 2^{-n-1}\delta.$ 
Since  $ \y \cE \yy = \infty$ the condensation property 
(cf. \cite[Lemma 3.1.4]{pfxm1}) yields  
the existence of a $K_0 \in \cE $ such that 
\[ \sum_{J \in G_n ( K_0 , \cE )} |J| \ge (1-\vare) |K_0| . \]
Examining the Gamlen-Gaudet construction \cite{gg} as (for example)
presented in \cite[Proposition 3.1.6]{pfxm1}, we obtain a family 
$(\cB_I)_{I\in\cD_n}$ of collections of dyadic intervals such that
\begin{enumerate}[(i)]
\item $\cB_I \subseteq K_0 \cap \cE$,
\item the elements of $\cB_I$ are pairwise disjoint,
\item for $B_I := \bigcup_{K \in \cB_I} K$ one has that $B_I\cap B_J=\emptyset$
      if and only if $I\cap J=\emptyset$, and
      $B_I\subseteq B_J$ if and only if $I\subseteq J$,
\item for 
      \[ k_I := \sum_{K\in \cB_I} h_K \]
      and $I,I^-,I^+\in\cD_n$ such that $I^-$ is the right half of $I$ and
      $I^+$ the left half of $I$, one has
      $B_{I^-}\subseteq \{ k_I = -1 \}$ and 
      $B_{I^+}\subseteq \{ k_I =  1 \}$,
\item for $0\le k \le n$ and  $|I|=2^{-k}$ one has 
      \[ \frac{|K_0|}{2^k} - 2 \vare |K_0| \le |B_I| \le  \frac{|K_0|}{2^k} .\]
\end{enumerate}
As a consequence 
\[ \frac{1-\delta}{2^n} |K_0| \le |B_I| \le  \frac{|K_0|}{2^n}
   \sptext{1}{for}{1}
   |I|=2^{-n} \]
and
\[  (1-\delta) |K_0| \le \sum_{|I|=2^{-n}} |B_I| \le |K_0|.\]
Choose measurable subsets $A_I \subseteq B_I$ for $|I|=2^{-n}$ such that
\begin{enumerate}[(a)]
\item $|A_I|=  (1-\delta) 2^{-n} |K_0|$,
\item the $k_I$ restricted to $A_I$ are symmetric.
\end{enumerate}
Let $S:=  \bigcup_{|I|=2^{-n}} A_I$, so that $|S|=(1-\delta)|K_0|$, and
\[ A_I := B_I \cap S 
   \sptext{1}{for all (remaining)}{1}
   I\in\cD_n. \]
When restricted to the probability  space $ (S,\frac{dt}{|S|}) $ 
the system  $(k_I)_{I\in\cD_n}$ has the same joint 
distribution as the usual Haar system $(h_I)_{I\in\cD_n}$ on the unit interval.
 Hence, as a consequence of the Gamlen-Gaudet construction, 
we obtain that 
\equa
      \left \| \sum_{I\in\cD_n} \frac{h_I}{|I|^{1/p}} x_I \right \|_{L_X^p}
& = & \left \| \sum_{I\in\cD_n} \frac{k_I}{|I|^{1/p}} x_I 
      \right \|_{L_X^p\left (S,\frac{dt}{|S|}\right ) } \\
& = & \left ( \frac{1}{|S|} \right )^\frac{1}{p}
      \left \| \sum_{I\in\cD_n} \frac{k_I}{|I|^{1/p}} x_I \right \|_{L_X^p(S,dt)} \\
&\le& \left ( \frac{1}{|S|} \right )^\frac{1}{p}
      \left \| \sum _{ I \in \cD _n}   
      \sum_{K\in \cB_I }  \frac{h_K}{|K|^{1/p} } 
      \left( \frac{|K|^{1/p}}{|I|^{1/p} } x_I \right) \right \|_{L_X^p}.
\tion
Recall that we selected the collection  $\cB_I$ as a sub-collection of $\cE$.
Using our hypothesis concerning $\cE $ and $X$, we obtain an upper estimate 
for the last term as follows,
\equa
      c \left( \frac{1}{|S|}\right )^\frac{1}{p}
      \left ( \sum _{ I \in \cD _n} \sum_{K\in \cB_I }   
      \frac{|K|}{|I|} \|x_I\|^p \right )^\frac{1}{p}
& = & c \left ( \sum _{ I \in \cD _n}     
      \frac{|B_I|}{|I| |S|} \|x_I\|^p  \right )^\frac{1}{p} \\
& = & c \left ( \sum _{ I \in \cD _n}     
      \frac{|B_I|}{|I| (1-\delta)|K_0|} \|x_I\|^p \right )^\frac{1}{p}  \\
&\le& \frac{c}{(1-\delta)^\frac{1}{p}}
      \left ( \sum _{ I \in \cD _n} \|x_I\|^p \right )^\frac{1}{p}.
\tion
Letting $\delta \downarrow 0$ yields our statement.
\endproof



\textbf{Addresses}
\parindent0em

Department of Mathematics and Statistics \\
P.O. Box 35 (MaD) \\
FIN-40014 University of Jyv\"askyl\"a \\
Finland 

\medskip
Department of Analysis\\
J. Kepler University\\
A-4040 Linz\\
Austria

\end{document}